\documentstyle[amssymb,amsfonts,12pt]{amsart}

\newcommand{\Chow}{\mbox{\it Chow}\,}
\newcommand{\Chowo}{{\mbox{\it Chow}}^\circ}

\newcommand{\DOT}{\setlength{\unitlength}{1pt}\begin{picture}(2.5,2)(1,1)
\put(1,2){\circle*{2}}\end{picture}}

\newcommand{\Edot}{{E_{\DOT}}}

\newcommand{\Fdot}{{F\!_{\DOT}}}
\newcommand{\Fpdot}{{{F\!_{\DOT}}'}}

\newcommand{\QED}{
\setlength{\unitlength}{1.0pt}%
\begin{picture}(7.5,7.5)
\put(0,-2.5){\rule{2.5pt}{5pt}}
\put(0,2.5){\rule{5pt}{2.5pt}}
\put(0,5){\rule{7.5pt}{2.5pt}}
\end{picture}}

\newcommand{\Span}[1]{\langle #1 \rangle}

\begin{document}

\title[real enumerative geometry and effective algebraic
equivalence]{Real Enumerative geometry and effective \\ algebraic equivalence}

\author{Frank Sottile}

\address{
        Department of Mathematics\\
        University of Toronto\\
        100 St. George Street\\
	Toronto, Ontario  M5S 3G3\\
	Canada\\
	(416) 978-4031}
\email{sottile@@math.toronto.edu}
\date{5 February 1996}
\thanks{Research supported in part by NSERC grant \# OGP0170279}

\maketitle

\section{Introduction}
Determining the common zeroes of a set of polynomials
is further complicated over non-algebraically closed fields such
as the real numbers.
A special case is whether a  problem of enumerative geometry can have all
its solutions be real.
We call such a problem {\em fully real}.

Little is known about enumerative geometry from this perspective.
A standard proof
of B\'ezout's Theorem shows the problem of intersecting hypersurfaces in
projective space is fully real.
Khovanskii~\cite{Khovanskii_fewnomials} considers
intersecting hypersurfaces in a torus defined by few monomials
and shows the real zeros are at most a fraction of the
complex zeroes.
Fulton, and more recently, Ronga, Tognoli and Vust~\cite{Ronga_Tognoli_Vust}
have shown the problem of 3264 plane conics tangent to five given
conics is fully real.
The author~\cite{sottile_real_lines} has shown
all problems of enumerating lines incident on linear
subspaces of projective space are fully real.

There are few methods for studying this phenomenon.
We ask:
How can the knowledge that one enumerative problem is fully
real be used to infer that a related problem is fully real?
We give several procedures to accomplish this inference and
examples of their application, lengthening the list of enumerative
problems known to be fully \smallskip real.

We study intersections of any dimension,
not just the zero dimensional intersections of enumerative problems.
Our technique is to deform general intersection
cycles into simpler cycles.
This modification of the classical method of degeneration
was used by Chiavacci and
Escamilla-Castillo~\cite{Chiavacci_Escamilla-Castillo}
to investigate these questions.

Let $\alpha_1,\ldots,\alpha_a$ be cycle classes spanning
the Chow ring of a smooth
variety.
For cycle classes $\beta_1,\ldots,\beta_b$,
there exist integers $c_i$ for $i = 1,\ldots, a$ such that
$$
\beta_1\cdots\beta_b \ =\
\sum_{i=1}^a c_i\cdot \alpha_i.
$$
When the $c_i\geq 0$,
this product formula has a geometric interpretation.
Suppose $Y_1,\ldots,Y_b$ are cycles representing the classes
$\beta_1,\ldots,\beta_b$ which meet generically transversally in a
cycle $Y$.
Then $Y$ is algebraically equivalent to
$Z:=Z_1\cup\cdots\cup Z_a$, where
$Z_i$ has $c_i$ components, each representing the cycle class $\alpha_i$.
This algebraic equivalence is effective if there is a family of cycles
containing both $Y$ and $Z$ whose general member is
a generically transverse intersection of cycles representing the
classes $\beta_1,\ldots,\beta_b$.
If the cycles $Y_1,\ldots,Y_b$ and each component of $Z$ are
defined over ${\Bbb R}$ and both $Y$ and $Z$ are in the same connected
component of the real points of that family,
then the effective algebraic equivalence is real.

Real effective algebraic equivalence
can be used to show an enumerative problem is fully real, or
more generally, to obtain lower bounds on the maximal
number of real solutions.
Suppose the cycles
$Y_1,\ldots,Y_b, W_1,\ldots, W_c$ give an enumerative problem
and the problem obtained by substituting $Z$ for $Y_1,\ldots,Y_b$
has at least $d$ real solutions.
Then there exist real cycles  $Y'_1,\ldots,Y'_b$ such that
the original problem (with $Y'_i$ in place of $Y_i$)
has at least $d$ real solutions,
since the number of real solutions is preserved by
small real deformations.

Sections 2 through 5 introduce and develop our basic notions and techniques.
Subsequent sections are devoted to elaborations and
applications of these ideas.

In Section 6, we  prove that any enumerative problem on a
flag variety involving five Schubert varieties, three of
which are special Schubert varieties, is fully real.
Given a map  $\pi: Y\rightarrow X$ with equidimensional fibres,
real effective algebraic equivalence on $X$ and $Y$ is compared
in Section 7 and used in Sections 8 and 9
to show that many Schubert-type enumerative problems
in two classes of flag varieties are fully real.

A proof of B\'ezout's Theorem in Section 10 suggests  another
method for obtaining fully real enumerative problems.
This is applied in Section 11 to show the enumerative problem of
$(n-2)$-planes in
${\Bbb P}^n$ meeting $2n-2$ rational normal curves is fully real.

The author thanks Bernd Sturmfels for encouraging these investigations.

\section{Intersection Problems}

\subsection{Conventions}
Varieties are reduced,  complex, and defined over  the
real numbers ${\Bbb R}$.
Let  $X$ and $Y$ denote smooth
projective varieties and $U$, $V$, and $W$ smooth quasi-projective
varieties.
Equip the real points $X({\Bbb R})$ of $X$
with the classical topology.
Let $A^*X$ be the Chow ring of cycles
modulo algebraic equivalence.

Two subvarieties meet {\em generically transversally}
if they meet  transversally along a dense subset of each
component of their intersection.
Such an intersection scheme is reduced at the generic point of each
component, or generically reduced.
A subvariety $\Xi\subset U\times X$ (or $\Xi\rightarrow U$)
with generically reduced equidimensional fibres
over a smooth base
$U$ is a family of {\em multiplicity free cycles on $X$ over $U$.}
All fibres of $\Xi$ over $U$ are algebraically equivalent, and we say
$\Xi\rightarrow U$ {\em represents} that algebraic equivalence class.

\subsection{Chow varieties}\label{sec:Chow}
Positive cycles on $X$ of a fixed dimension and degree are parameterized
by the Chow variety  of $X$.
We  suppress the dependence on dimension and degree
and write $\Chow X$ for any Chow variety of $X$.
The open Chow variety $\Chowo X$ is the open subset of $\Chow X$
parameterizing multiplicity free cycles on $X$.
There is a tautological family
$\Phi \rightarrow \Chowo X$
of cycles on $X$ with the property
that $\zeta\in\Chowo X$ represents the the fundamental cycle of
the fibre $\Phi_\zeta$.

Let $\Xi\rightarrow U$ be a family of multiplicity
free cycles on $X$.
The association of a point $u$ of $U$ to the fundamental cycle of
the fibre $\Xi_u$ defines the
{\em fibre function} $\phi$,  which is algebraic
on a dense open subset $U'$ of $U$.
If $U$ is a  curve, then $U=U'$.

\subsection{Proposition}
{\em $\phi(U')$ is dense in the set $\phi(U)$.}\medskip

\noindent{\bf Proof:}
Let $u\in U$ and $C\subset U$ be a smooth curve with $u\in C$ and
$C-\{u\}\subset U'$.
Such a curve is not necessarily closed in
$U$, but is the smooth points of a closed subvariety.
The fibre function $\phi|_C$ of $\Xi|_C \rightarrow C$
is algebraic, hence
$\phi(u)\in \overline{\phi(C-\{u\})} \subset \overline{\phi(U')}$.
\QED\vspace{10pt}

Two families  $\Xi\rightarrow  U$ and $\Psi\rightarrow V$
of multiplicity free cycles on $X$ are {\em equivalent} if
$\overline{\phi(U)} = \overline{\phi(V)}$, that is, if
they have essentially the same cycles.
Our results remain valid when
one family of cycles is replaced by an equivalent family,
perhaps with the additional assumption that
$\overline{\phi(U({\Bbb R}))} = \overline{\phi(V({\Bbb R}))}$.

The varieties $\Chow X$ and $\Chowo X$ as well as $U'$ and the
morphism $\phi: U' \rightarrow \Chow X$ are defined over
${\Bbb R}$~(\cite{Samuel}, \S I.9).
We use $\phi$ to denote all fibre functions.
Any ambiguity may be resolved by context.

\subsection{Intersection Problems}\label{sec:intersection_problems}
For  $1\leq i\leq b$,
let $\Xi_i\rightarrow U_i$  be a family of multiplicity free cycles
on $X$.
Let $U\subset \prod_{i=1}^b U_i$ be the locus where the fibres of the
product family $\prod_{i=1}^b \Xi_i$ meet the (small) diagonal
$\Delta^b_X$ of $X^b$ generically transversally.
Equivalently, $U$ is the locus where fibres of
$\Xi_1,\ldots,\Xi_b$ meet generically transversally in $X$.
If $U$ is nonempty, then
$\Xi_1,\ldots,\Xi_b$ give a {\em (well-posed) intersection problem}.

Given an intersection problem as above, let
$\delta:X\stackrel{\sim}{\longrightarrow}\Delta^b_X\subset X^b$ and set
$\Xi$ to be
$$
\Xi \ :=\ (1_U\times \delta)^* \prod_{i=1}^b \Xi_i \ \subset\ U\times X,
$$
a family of multiplicity
free cycles on $X$ over $U$.
We often suppress the dependence on the original families and write
$\Xi\rightarrow U$ for this intersection problem.

Not all collections of families of cycles give  well-posed
intersection problems,
some transversality is needed to guarantee $U$ is nonempty.
When a reductive group acts transitively on $X$, Kleiman's
Transversality Theorem~\cite{Kleiman} has the following consequence.

\subsection{Proposition}\label{prop:transitive_action}
{\em
Suppose a reductive group acts transitively on $X$, $\Xi_1$ is a constant
family,  and for
$2\leq i\leq b$, $\Xi_i$ is equivalent to a family of multiplicity free cycles
stable under that action.
Then $\Xi_1,\ldots,\Xi_b$ give  a well-posed intersection problem.
}\medskip

Grassmannians and flag varieties have such an action.
For these, we suppose all families of cycles
are stable under that action, and thus
give well-posed intersection problems.

Suppose a reductive group acts on $X$ with a dense open orbit $X'$.
For example, if $X$ is a toric variety, or more generally, a
spherical variety~\cite{Brion_spherical_introduction,%
Knop_spherical_expository,Luna_Vust_Plongements}.
Each family may be stable under that action, but the collection need not
give a well-posed intersection problem as
Kleiman's Theorem~\cite{Kleiman} only guarantees transversality in $X'$.
However, it is often the case that only points of
intersection in $X'$ are desired,
and suitable blow up of $X$ or a different equivariant compactification
of $X'$ exists on which
the corresponding intersection problem is well-posed
(see~\cite{Fulton_introduction_intersection}, \S 1.4
or~\cite{Fulton_intersection}, \S 9 and \S 10.4).

\section{effective algebraic equivalence}
\label{sec:effective_algebraic_equivallences}

Let $\alpha_1,\ldots,\alpha_a$ be distinct
additive generators of $A^*X$,
and for $1\leq i\leq a$, suppose
$\Psi(\alpha_i)\rightarrow V(\alpha_i)$ is a family of multiplicity free
cycles on $X$ representing the cycle class $\alpha_i$.
When $X$ is a Grassmannian or flag variety,
$\alpha_1,\ldots,\alpha_a$ will be the Schubert classes,
and $\Psi(\alpha_i)\rightarrow V(\alpha_i)$ the corresponding families
of Schubert varieties.

A family of multiplicity free cycles
$\Xi\subset U\times X$ has an {\em effective algebraic equivalence}
with {\em witness} $Z\in \overline{\phi(U)}\cap \Chowo X$
if each (necessarily multiplicity free) component of $Z$ is a fibre
of some  family $\Psi(\alpha_i)$.
This  effective algebraic equivalence is {\em real}  if
$Z\in \overline{\phi(U({\Bbb R}))}$  and
each component of $Z$ is a fibre over a real point of some $V(\alpha_i)$.
An intersection problem
$\Xi_1,\ldots,\Xi_b$ has
{\em  (real) effective algebraic equivalences}
if its  family of intersection cycles $\Xi\rightarrow U$
has (real) effective algebraic equivalences.

\subsection{Products in $A^*X$}
\label{sec:products}
Suppose $\beta_1,\ldots,\beta_b$ are classes from
$\{\alpha_1,\ldots,\alpha_a\}$
and the families  $\Psi(\beta_1),\ldots,\Psi(\beta_b)$
give an intersection problem $\Psi\rightarrow V$.
We say  $\Psi\rightarrow V$ is given by $\beta_1,\ldots,\beta_b$.
Suppose $\Psi\rightarrow V$  has an
effective algebraic equivalence with witness $Z$.
Fibres of
$\Psi\rightarrow V$ are generically transverse intersections
of fibres of $\Psi(\beta_1),\ldots,\Psi(\beta_b)$,  and so have
cycle class $\beta_1\cdots\beta_b$.
As $Z\in \overline{\phi(V)}$, this equals the cycle class of $Z$,
which is $\sum_{i=1}^a c_i \alpha_i$,
where $c_i$ counts the components of $Z$
lying in the family $\Psi(\alpha_i)$.
Thus in $A^*X$, we have
$$
\hspace{2.5in}
\beta_1\cdots\beta_b \ =\ \sum_{i=1}^a c_i \alpha_i.
\hspace{2.2in}
(\ref{sec:products})
$$

To compute products in $A^*X$, classical geometers would try to
understand a generically transverse intersection of degenerate
cycles in special position, as a generic intersection cycle is
typically too difficult to describe.
Effective algebraic equivalence extends this method of degeneration
by also considering
limiting positions of such intersection cycles as the subvarieties
degenerate further, attaining excess intersection.

\subsection{Pieri-type intersection problems}\label{sec:pieri_type}
A Schubert subvariety of a flag variety is determined by a
complete flag $\Fdot$ and a
coset $w$ of a parabolic subgroup in the symmetric group.
Thus Schubert classes $\sigma_w$ are indexed by these cosets
and families $\Psi_w$ of Schubert varieties have base
${\Bbb F}\ell$, the variety of complete flags.

A {\em special Schubert subvariety} of a Grassmannian
is the locus of planes meeting a fixed linear subspace non-trivially,
or the image of such a subvariety in the dual Grassmannian.
More generally, a special Schubert subvariety of a flag variety
is the pullback of a special Schubert subvariety from
a Grassmannian projection.
If $m$ is the index of a special Schubert class, then the Pieri-type
formulas
of~\cite{Lascoux_Schutzenberger_polynomes_schubert,Sottile_Pieri_Schubert}
show that for any $w$,
there exists a subset $I_{m,w}$ of these cosets
such that
$$
\hspace{2.5in}
\sigma_m\cdot \sigma_w\ =\ \sum_{v\in I_{m,w}}\sigma_v.
\hspace{2.2in}
(\ref{sec:pieri_type})
$$

\subsection{Theorem}\label{thm:pieri_effective_equivalences}
{\em
The intersection problem $\Xi\rightarrow U$ given by the classes
$\sigma_m$ and $\sigma_w$
has real effective algebraic equivalences.
}\medskip

\noindent{\bf Proof:}
The Borel subgroup $B$ of $GL_n{\Bbb C}$ stabilizing a real
complete flag $\Fdot$
acts on the Chow variety with fixed points the
$B$-stable cycles, which are sums of
Schubert varieties determined by $\Fdot$.
As Hirschowitz~\cite{Hirschowitz} observed,
$\overline{\phi(U)}$ is $B$-stable,  and must
contain a fixed point (\cite{Borel_groups}, III.10.4).
In fact,  if $\Fpdot$ is a real flag in linear general position with $\Fdot$,
then the $B({\Bbb R})$-orbit of $\Omega_m\Fdot\bigcap \Omega_w\Fpdot$
is a subset of $\phi(U({\Bbb R}))$.
Moreover its closure has a $B({\Bbb R})$-fixed point, as the
proof in~\cite{Borel_groups}
may be adapted to show that complete $B({\Bbb R})$-stable
real analytic sets have fixed points.
Since the coefficients of the sum  (\ref{sec:pieri_type}) are
all 1, $\sum_{v\in I_{m,w}}\Omega_w\Fdot$ is the only $B({\Bbb R})$-stable
cycle in its algebraic equivalence class, and therefore
$$
\sum_{v\in I_{b,w}}\Omega_w\Fdot\ \in\ \overline{\phi(U({\Bbb R}))}.
\qquad\qquad\qquad\qquad\QED
$$

\section{Fully real enumerative problems}

An {\em enumerative problem} of {\em degree} $d$ is an
intersection problem $\Xi\rightarrow U$
with zero-dimensional fibres of cardinality $d$.
An enumerative problem is {\em fully real} if there exists
$u\in U({\Bbb R})$ with
all points in the fibre $\Xi_u$ real.
In this case, $u = (u_1,\ldots,u_b)$ with $u_i\in U_i({\Bbb R})$
and the cycles $(\Xi_1)_{u_1},\ldots, (\Xi_b)_{u_b}$  meet transversally with
all points of intersection real.

\subsection{Theorem}\label{thm:real_closure}{\em
An enumerative problem $\Xi\rightarrow U$ is fully real if and
only if it has real effective algebraic equivalences.
That is, if and only if there exists a
point $\zeta\in \overline{\phi(U({\Bbb R}))}$
representing distinct real points.
}\medskip

\noindent{\bf Proof:}
The forward implication is a consequence of the definition.
For the reverse,
let $d$ be the degree of $\Xi\rightarrow U$.
Then $\phi:U\rightarrow S^dX$, the Chow variety of effective degree $d$
zero cycles on $X$.
The real points $S^dX({\Bbb R})$ of the Chow variety
represent degree $d$ zero cycles stable
under complex conjugation.
Its dense set of multiplicity free cycles
have an open subset ${\cal M}$ parameterizing cycles of
distinct real points, and $\zeta\in {\cal M}$.
Thus $\phi(U({\Bbb R}))\cap {\cal M}\neq \emptyset$,
which implies $\Xi\rightarrow U$ is  fully real.
\QED\vspace{10pt}

The set of witnesses to $\Xi\rightarrow U$ being fully real
contains an open subset $\phi^{-1}({\cal M})\bigcap U({\Bbb R})$.

\section{Curve selection}

Subsequent sections use real effective algebraic
equivalence for one or more families to infer results
about related  families.
While intuition supports the claim that the
functions we define between Chow varieties are algebraic
(or at least continuous),
we are unaware of general results verifying this intuition.
An obvious obstruction is that the Chow variety does not represent a functor.
However, weaker claims suffice.
Our tool is the Curve Selection
Lemma~\cite{Benedetti_Risler} of real
semi-algebraic geometry, in the following guise:

\subsection{Curve Selection Lemma}\label{lemma:curve_selection}
{\em Let $V$ be a real variety and $R\subset V({\Bbb R})$ a
semi-algebraic subset.
If $\zeta \in \overline{R}$, then there is a real algebraic map
$f: C\rightarrow V$ with   $C$ a smooth curve,
and a point $s$ on a connected arc $S$ of
$C({\Bbb R})$ such that $f(S -\{s\}) \subset R$
and $f(s) = \zeta$.\medskip }

\noindent{\bf Proof:}
By the Curve Selection Lemma
(\cite{Benedetti_Risler}, 2.6.20),
there exists a semi-algebraic function
$g:[0,1]\rightarrow \overline{R}$ with $g(0)=\zeta$,
$g(0,1]\subset R$, and $g$ a real analytic homeomorphism onto its image
in $\overline{R}$.
Let $C^\circ$ be the Zariski closure of $g[0,1]$ in $V$,
and $f:C\rightarrow C^\circ$ its normalization.
Let $S\subset C({\Bbb R})$ be a connected arc of $f^{-1}(g[0,1])$ whose image
contains $g(0)$ and let $s \in S\cap f^{-1}(g(0))$.
\QED

\section{Pieri-type enumerative problems}

\subsection{Theorem}\label{thm:pieri_schubert}
{\em  Any enumerative problem in any flag variety involving
five Schubert varieties, three of which are special, is fully real.}
\medskip

This generalizes Theorem 5.2 of~\cite{sottile_explicit_pieri},
the analogous result for Grassmannians.
It requires an additional transversality result.

\subsection{Lemma}\label{lemma:special_transversality}
{\em
Let $(w_1,w_2)$ and $(v_1,v_2)$ be indices of Schubert subvarieties of
a flag variety, with $w_1$ and $w_2$ (respectively $v_1$ and $v_2$)
defining defining  Schubert varieties of
the same dimension.
Suppose $m$ is the index of a special Schubert subvariety
such that $(w_1,v_1,m)$ gives an enumerative problem.
If $(w_1,v_1)\neq (w_2,v_2)$, and $\Fdot, \Fpdot$
are complete flags
in linear general position, then
there is an open set $V$ of the variety
${\Bbb F}\ell$ of
complete flags consisting of flags $\Edot$ such that
$$
\Omega_{w_1}\Fdot \bigcap \Omega_{v_1}\Fpdot \bigcap \Omega_m \Edot
\qquad\mbox{and}\qquad
\Omega_{w_2}\Fdot \bigcap \Omega_{v_2}\Fpdot \bigcap \Omega_m \Edot
$$
are transverse intersections which coincide only when empty.
}\medskip

If the three flags are real, then  a nonempty intersection
as above is a single real point.
\medskip

\noindent{\bf Proof:}
By Kleiman's Theorem ~\cite{Kleiman},
there is an open subset $U$ of
${\Bbb F}\ell\times{\Bbb F}\ell\times{\Bbb F}\ell$
consisting of triples $(\Fdot,\Fpdot,\Edot)$ such that each
intersection is transverse and so
is either empty or a single point,
by the Pieri-type formulas
of~\cite{Lascoux_Schutzenberger_polynomes_schubert,Sottile_Pieri_Schubert}.
Suppose neither is empty.

Similarly, there is an open subset $V$ of triples for which
$$
\left(\Omega_{w_1}\Fdot \bigcap \Omega_{w_2}\Fdot \right)
 \bigcap \left(\Omega_{v_1}\Fpdot\bigcap \Omega_{v_2}\Fpdot\right)
 \bigcap \Omega_b \Edot
$$
is proper.
Since $(w_1,v_1)\neq (w_2,v_2)$, it is empty.
Thus for triples $(\Fdot,\Fpdot,\Edot)\in U\cap V$,
$$
\Omega_{w_1}\Fdot \bigcap \Omega_{v_1}\Fpdot \bigcap \Omega_b \Edot\  \neq\
\Omega_{w_2}\Fdot \bigcap \Omega_{v_2}\Fpdot \bigcap \Omega_b \Edot.
$$

The lemma follows, as $U\cap V$ is stable under the
diagonal action of $GL_n{\Bbb C}$
and the set of pairs $(\Fdot,\Fpdot)$ in linear
general position is the open $GL_n{\Bbb C}$-orbit
in ${\Bbb F}\ell\times{\Bbb F}\ell$.
\QED

\subsection{Proof of Theorem~\ref{thm:pieri_schubert}:}
Let $\Xi_1,\Xi_2,\Xi_3, \Gamma_1$, and $\Gamma_2$ be families of
Schubert varieties representing the classes
$\sigma_{m_1}, \sigma_{m_2}, \sigma_{m_3}, \sigma_{w_1}$, and $\sigma_{w_2}$.
Suppose $\sigma_{m_1}, \sigma_{m_2}$, and $\sigma_{m_3}$ are special Schubert
classes, and these families give an enumerative problem
$\Xi\rightarrow U$.

By \S \ref{sec:pieri_type}, for each $i=1,2$, the intersection problem
$\Psi_i\rightarrow V_i$ given by the families $\Xi_i$ and $\Gamma_i$
has a real effective algebraic equivalence
with witness
$\sum_{v_i\in I_{m_i,w_i}} \Omega_{v_i}\Fdot$, for any real flag $\Fdot$.
Let  $\Fdot$ and $\Fpdot$ be real flags in linear general position
and set
$$
Z_1 \ := \sum_{v_1\in I_{m_1,w_1}} \Omega_{v_1} \Fdot
\qquad\mbox{and}\qquad
Z_2 \ := \sum_{v_2\in I_{m_2,w_2}} \Omega_{v_2} \Fpdot.
$$

For $i=1,2$,
let $\phi_i$ be the fibre function for $\Psi_i\rightarrow V_i$.
Then $Z_i\in \overline{\phi_i(V_i({\Bbb R}))}$
and by  Lemma~\ref{lemma:curve_selection}, there is a
map $f_i : C_i \rightarrow \overline{\phi_i(V_i)}$
with $C_i$ a smooth curve, and a point $s_i$
on a connected arc $S_i$ of
$C_i({\Bbb R})$ such that $f(S_i -\{s_i\})\subset\phi_i(V_i({\Bbb R}))$
and $f_i(s_i) = Z_i$.
Then $f_i^*\Phi\rightarrow C_i$ is a family of multiplicity free cycles,
where $\Phi\rightarrow \Chowo X$ is the tautological family.

Considering pairs of components of $Z_1$ and $Z_2$ separately,
Lemma~\ref{lemma:special_transversality} shows there is a real flag
$\Edot$ such that $Z_1\bigcap Z_2 \bigcap \Omega_{m_3} \Edot$
is a transverse intersection all of whose points are real.
Thus $f_1^*\Phi \rightarrow C_1$,  $f_2^*\Phi \rightarrow C_2$, and
$\Xi_3\rightarrow {\Bbb F}\ell$ give a well-posed fully real enumerative
problem $\Psi \rightarrow V$, as
$(s_1,s_2,\Edot) \in V({\Bbb R})$.

Let ${\cal M}$ be the open subset of the real points of the
Chow variety parameterizing cycles consisting entirely of real points.
Then $\phi(s_1,s_2,\Edot)\in {\cal M}$
and so $\phi^{-1}({\cal M})$ meets
$R:= (S_1-\{s_1\})\times(S_2-\{s_2\})\times \{\Edot\}$.
However, fibres of $\Psi$ over points of
$R$ are fibres of
$\Xi$ over points of $U({\Bbb R})$,
showing $\Xi\rightarrow U$ to be fully real.
\QED

\section{Fibrations}
Suppose $\pi: Y\rightarrow X$ has equidimensional fibres.
If $\Xi\rightarrow  U$ is a family of multiplicity free cycles
on $X$ representing the cycle class $\alpha$, its pullback
$\pi^*\Xi := (\pi\times 1_U)^{-1}\Xi\rightarrow U$ is a family of
multiplicity free cycles on
$Y$ representing the cycle class $\pi^*\alpha$.

Suppose $\alpha_1\ldots,\alpha_a$  generate $A^*X$ additively
and $\Psi(\alpha_1),\ldots,\Psi(\alpha_a)$ are
families of cycles representing these generators.
The classes $\pi^*\alpha_1,\ldots,\pi^*\alpha_a$ generate the image of
$A^*X$ in $A^*Y$  and are represented by the families
$\pi^*\Psi(\alpha_1),\ldots,\pi^*\Psi(\alpha_a)$.
Effective algebraic equivalence is preserved by pullbacks:

\subsection{Theorem}\label{thm:pullbacks}
{\em
If $\,\Xi\rightarrow U$ is a family of multiplicity free cycles on
$X$ having effective algebraic equivalences with witness $Z$, then
$\pi^*\Xi\rightarrow U$ is a family of multiplicity free cycles on
$Y$ having effective algebraic equivalences with witness $\pi^{-1}Z$.
Likewise, if $\,\Xi\rightarrow U$ has real effective
algebraic equivalences, then so does $\pi^*\Xi\rightarrow U$.
}\medskip

Associating a cycle $Z$ on $X$ to $\pi^{-1}Z\subset Y$
defines a function $\pi^*: \Chow X \rightarrow \Chow Y$.
If $\phi$ is the fibre function of $\Xi\rightarrow U$, then
$\pi^*\circ \phi$ is the fibre function of $\pi^*\Xi\rightarrow U$.
Letting $W=\phi(U')$ and $R = \phi(U'({\Bbb R}))$, we see that
Theorem~\ref{thm:pullbacks} is a consequence of the following
lemma.

\subsection{Lemma}\label{lemma:pullback}
{\em
Let $W\subset \Chowo X$ be constructible and $V := \overline{W}\cap \Chowo X$.
Then $\pi^*(V)\subset \overline{\pi^*(W)}$ in $\Chow Y$.
Likewise, if $R\subset \Chowo X({\Bbb R})$ is semi-algebraic
and $Q:=\overline{R}\cap \Chowo X({\Bbb R})$, then
$\pi^*(Q)\subset \overline{\pi^*(R)}$ in $\Chow Y({\Bbb R})$.
}\medskip

\noindent{\bf Proof:}
For the first part,  let $\zeta\in V$.
We show $\pi^*(\zeta) \in \overline{\pi^*(W)}$.

Let  $C^\circ\subset \Chowo X$ be an irreducible curve with
$\zeta\in C^\circ$ and $C^\circ-\{\zeta\} \subset W$.
Let $f: C \rightarrow C^\circ$ be its normalization and let
$s\in f^{-1}(\zeta)$.
Let $\Phi \subset \Chowo X \times X$ be the tautological family.
Then $f^*\Phi$ is a family of multiplicity free cycles on $X$
with fibre function $f$.
Similarly, $\pi^*\circ f$ is the fibre function of the family
$\pi^*(f^*\Phi)$ of multiplicity free cycles on $Y$ over the smooth curve $C$.
As noted in \S \ref{sec:Chow},  this implies $\pi^*\circ f$ is algebraic,
and so
$\pi^*(\zeta) \in \pi^*(f(C)) \subset \overline{\pi^*(W)}$, since
$\pi^*(f(C-\{f^{-1}(\zeta)\}))\subset \pi^*(W)$.

For the second part, suppose $R\subset \Chowo X({\Bbb R})$
and $\zeta \in Q = \overline{R}\bigcap\Chowo X$.
By Lemma~\ref{lemma:curve_selection}, there is a smooth curve $C$,
a connected arc $S\subset C({\Bbb R})$, a point $s\in S$, and
an algebraic map $f: C\rightarrow \Chowo X$ such that
$f(s) = \zeta$ and $f(S-\{s\}) \subset R$.
Arguing as above shows
$\pi^*(\zeta) \in \pi^*(f(S)) \subset \overline{\pi^*(R)}$.
\QED

\section{Schubert-type enumerative problems in
${\Bbb F}\ell_{0,1}{\Bbb P}^n$ are fully real}

The variety ${\Bbb F}\ell_{0,1}{\Bbb P}^n$ of partial
flags $q\in l \subset {\Bbb P}^n$ with
$q$ a point  and $l$ a line has  projections
$$
p\ :\ {\Bbb F}\ell_{0,1}{\Bbb P}^n \longrightarrow {\Bbb P}^n
\qquad\mbox{and}\qquad
\pi\ :\ {\Bbb F}\ell_{0,1}{\Bbb P}^n\longrightarrow
{\Bbb G}_1{\Bbb P}^n,
$$
where ${\Bbb G}_1{\Bbb P}^n$ is  the Grassmannian of lines
in ${\Bbb P}^n$.

A Schubert subvariety of ${\Bbb G}_1{\Bbb P}^n$ is determined by
a partial flag $F\subset P$ of ${\Bbb P}^n$:
$$
\Omega(F,P)\ :=\ \{l \in {\Bbb G}_1{\Bbb P}^n\,|\,
l\cap F \neq \emptyset\ \mbox{\ and }\ l\subset P\}.
$$
If $F$ is a hyperplane of $P$, then $\Omega(F,P) = {\Bbb G}_1P$, the
Grassmannian of lines  in $P$.

In addition to $\pi^{-1}\Omega(F,P)$, there is another Schubert subvariety of
${\Bbb F}\ell_{0,1}{\Bbb P}^n$ which projects onto $\Omega(F,P)$ in
${\Bbb G}_1{\Bbb P}^n$:
$$
\widehat{\Omega}(F,P)\ :=\
\{(q,l) \in {\Bbb F}\ell_{0,1}{\Bbb P}^n\,|\,
q\in F \ \mbox{\ and }\ l\subset P\}.
$$

Any Schubert subvariety of ${\Bbb F}\ell_{0,1}{\Bbb P}^n$ is one of
$\Omega(F,P)$ or $\widehat{\Omega}(F,P)$, for suitable $F\subset P$.
The varieties $\widehat{\Omega}(F,P)$ have another description,
which is straightforward to verify:

\subsection{Lemma}\label{lemma:no_hats}
{\em
Let $N,P$ be subspaces of $\,{\Bbb P}^n$.
Then
$$
p^{-1}N \bigcap \pi^{-1}{\Bbb G}_1 P \ =\
\widehat{\Omega}(N\cap P,\,P),
$$
and, if $N$ and $P$ meet properly, this intersection
is generically transverse.
}

\subsection{Corollary}\label{cor:reduction}
{\em
Any Schubert-type enumerative problem on ${\Bbb F}\ell_{0,1}{\Bbb P}^n$
is equivalent to one involving only pullbacks of Schubert subvarieties
of $\,{\Bbb P}^n$ and ${\Bbb G}_1{\Bbb P}^n$.
}\medskip

The next lemma,  an exercise in linear algebra,
describes Poincar\'e duality for
Schubert subvarieties of ${\Bbb F}\ell_{0,1}{\Bbb P}^n$.

\subsection{Lemma}\label{lemma:poincare_duality}
{\em
Suppose a linear subspace $N$ meets a partial flag $F\subset P$
properly in ${\Bbb P}^n$.
If $\pi^{-1}\Omega(F,P)$ and $p^{-1}N$ have complimentary
dimension in ${\Bbb F}\ell_{0,1}{\Bbb P}^n$, then their intersection
is empty unless $F$ and $N\cap P$ are points.
In that case, they meet transversally in a single point and
$\pi^{-1}\Omega(F,P)\bigcap p^{-1}N
= (N\cap P,\ \Span{F, N\cap P})$.}

\subsection{Theorem}\label{thm:schubert_fully_real}
{\em
Any Schubert-type enumerative problem in ${\Bbb F}\ell_{0,1}{\Bbb P}^n$
is fully real.
}\medskip

\noindent{\bf Proof:}
By Corollary~\ref{cor:reduction}, it suffices to consider
enumerative problems involving only pullbacks of
Schubert subvarieties of ${\Bbb P}^n$ and ${\Bbb G}_1{\Bbb P}^n$.
Since the intersection of linear subspaces in ${\Bbb P}^n$
is another linear subspace, we may further suppose the enumerative
problem $\Xi\rightarrow U$ is given by families
$p^*\Xi_1,\pi^*\Xi_2,\ldots,\pi^*\Xi_b$, where $\Xi_1$ is the family of
subspaces of a fixed dimension in  ${\Bbb P}^n$ and
$\Xi_2,\ldots,\Xi_b$ are families of Schubert subvarieties
of ${\Bbb G}_1{\Bbb P}^n$.

By Theorem~C$'$ of~\cite{sottile_real_lines}, the intersection problem
$\Psi \rightarrow V$ on
${\Bbb G}_1{\Bbb P}^n$ given by $\Xi_2,\ldots,\Xi_b$ has
real effective algebraic equivalences.
Let $Z$ be a witness.
By Theorem~\ref{thm:pullbacks},
$\pi^*\Psi \rightarrow V$ has a real
effective algebraic equivalence with witness $\pi^* Z$.

By Lemma~\ref{lemma:curve_selection},
there is a real algebraic map
$f:C \rightarrow \overline{\phi(W)}\cap \Chowo {\Bbb F}\ell_{0,1}{\Bbb P}^n$
with $C$ a smooth curve, and a point $s$ on
a connected arc $S$ of $ C({\Bbb R})$ such that $f(s) = \pi^{-1}Z$ and
$f(S-\{s\})\subset \phi(V({\Bbb R}))$.
Let $\Phi \rightarrow \Chowo {\Bbb F}\ell_{0,1}{\Bbb P}^n$ be the tautological
family
and consider the family $f^*\Phi \rightarrow C$.
The fibre over $s$ of $f^*\Phi$ is $\pi^{-1}Z$.

Let ${\cal L}$ be the lattice of subspaces of ${\Bbb P}^n$ generated by
the (necessarily real) subspaces defining components of $Z$,
and let $N$ be a real subspace from the family $\Xi_1$ meeting all
subspaces of ${\cal L}$ properly.
By  Lemma~\ref{lemma:poincare_duality},
$p^{-1}N\bigcap \pi^{-1}Z$ is transverse with
all points of intersection real.
Thus there is a Zariski open subset $C'$ of $C$ such that fibres of
$f^*\Phi$ over $C'$ meet $p^{-1} N$ transversally.
Then $s\in \overline{(S-\{s\})\bigcap C'({\Bbb R})}$, so there is
a point  $t\in S-\{s\}$ such that
$p^{-1}N\bigcap(f^*\Phi)_t$ is transverse and consists
entirely of real points.
But $f(t)\in \phi(V({\Bbb R}))$, so $(f^*\Phi)_t=\Phi_{f(t)}$ is a fibre of
$\pi^*\Psi$ over $V({\Bbb R})$, and hence a
generically transverse intersection of real Schubert varieties
from the families $\pi^*\Xi_2,\ldots,\pi^*\Xi_b$.
Thus $\Xi\rightarrow U$ is fully real.
\QED

\subsection{Effective algebraic equivalence for
${\Bbb F}\ell_{0,1}{\Bbb P}^n$}
Any Schubert-type intersection problem on
${\Bbb F}\ell_{0,1}{\Bbb P}^n$ has real effective algebraic
equivalences.
We give an outline, as a complete analysis is lengthy and involves
no new ideas beyond~\cite{sottile_real_lines}.

By Corollary~\ref{cor:reduction}, it suffices
to consider intersection problems $\Xi\rightarrow U$ given by families
$p^*\Xi_1,\pi^*\Xi_2,\ldots,\pi^*\Xi_b$, where $\Xi_1$ is a family of
subspaces of a fixed dimension in  ${\Bbb P}^n$ and
$\Xi_2,\ldots,\Xi_b$ are families of Schubert subvarieties
of ${\Bbb G}_1{\Bbb P}^n$.

In~\cite{sottile_real_lines}, the intersection problem given by
$\Xi_2,\ldots,\Xi_b$ is shown to have real effective algebraic
equivalences with witness $Z$.
Let $\Psi \rightarrow V$ be the intersection problem given by
$p^*\Xi_1$ and the constant family $\pi^{-1}Z$.
Using Theorem~\ref{thm:pullbacks} and Lemma~\ref{lemma:curve_selection}
one may show
$$
\phi(V)\subset \overline{\phi(U)}
\qquad\mbox{and}\qquad
\phi(V({\Bbb R}))\subset \overline{\phi(U({\Bbb R}))}.
$$
It suffices to show  $\Psi\rightarrow V$ has real
effective algebraic equivalences.

A proof that $\Psi\rightarrow V$ has real effective algebraic
equivalences mimics the proof of Theorem~E of~\cite{sottile_real_lines},
with the following Lemma playing the role of
Lemma~2.4 of~\cite{sottile_real_lines}.

\subsection{Lemma}\label{lemma:reducible}
{\em Let $F,P,N$, and $H$ be linear subspaces of ${\Bbb P}^n$
and suppose that $H$ is a hyperplane containing neither  $P$
nor $N$, $F$ is a proper subspace of $P\cap H$, and $N$ meets $F$,
and hence $P$ properly.
Set $L = N\cap H$.
Then $\pi^{-1}\Omega(F,P)$ and $p^{-1}L$ meet
generically transversally,
$$
\pi^{-1}\Omega(F,P)\bigcap p^{-1}L\ =\
\widehat{\Omega}(N\cap F,\,P) \ +\
\pi^{-1}\Omega(F,P\cap H)\bigcap p^{-1}N,
$$
and the second term is itself an irreducible generically
transverse intersection.}
\medskip

The proof of this statement is almost identical to the proof of
Lemma~2.4 of~\cite{sottile_real_lines}.

\section{Some Schubert-type  enumerative problems in
${\Bbb F}\ell_{1,n-2}{\Bbb P}^n$}

The variety ${\Bbb F}\ell_{1,n-2}{\Bbb P}^n$ of partial flags
$l \subset \Lambda \subset {\Bbb P}^n$, where
$l$ is a line and $\Lambda$ an $(n-2)$-plane
has projections
$$
\pi: {\Bbb F}\ell_{1,n-2}{\Bbb P}^n \rightarrow {\Bbb G}_1{\Bbb P}^n
\qquad\mbox{and}\qquad
p: {\Bbb F}\ell_{1,n-2}{\Bbb P}^n\rightarrow {\Bbb G}_{n-2}{\Bbb P}^n,
$$
where ${\Bbb G}_{n-2}{\Bbb P}^n$ is the Grassmannian of $(n-2)$-planes in
${\Bbb P}^n$.

\subsection{Theorem}\label{thm:simple_schubert}
{\em
Any enumerative problem in ${\Bbb F}\ell_{1,n-2}{\Bbb P}^n$
given by pullbacks of Schubert subvarieties of $\,{\Bbb G}_1{\Bbb P}^n$
and ${\Bbb G}_{n-2}{\Bbb P}^n$
is fully real.
}\medskip

\noindent{\bf Proof:}
Suppose $\pi^*\Xi_1,\ldots,\pi^*\Xi_b, p^*\Gamma_1,\ldots,p^*\Gamma_c$
give an enumerative problem on ${\Bbb F}\ell_{1,n-2}{\Bbb P}^n$
where, for $1\leq i\leq b$, $\Xi_i$ is a family of Schubert subvarieties
of  ${\Bbb G}_1{\Bbb P}^n$
and for $1\leq j\leq c$, $\Gamma_i$ is a family of Schubert subvarieties
of  ${\Bbb G}_{n-2}{\Bbb P}^n$.

By Theorem~\ref{thm:pullbacks} and~\cite{sottile_real_lines},
$\pi^*\Xi_1,\ldots,\pi^*\Xi_b$ give an intersection
problem $\Psi_1\rightarrow V_1$
which has a real algebraic
equivalence with witness $Z_1$.
Identifying ${\Bbb P}^n$ with its dual projective space
gives an isomorphism
${\Bbb G}_{n-2}{\Bbb P}^n
\stackrel{\sim}{\rightarrow} {\Bbb G}_1{\Bbb P}^n$,
mapping Schubert subvarieties to Schubert subvarieties.
It follows that  $p^*\Gamma_1,\ldots,p^*\Gamma_c$ give an intersection problem
$\Psi_2\rightarrow V_2$ which has a real algebraic
equivalence with witness $Z_2$.
It suffices to show the enumerative problem
$\Psi\rightarrow V$ given by $\Psi_1$ and $\Psi_2$ is fully real.

Since $Z_1$ and $Z_2$ may be replaced by any translate
by elements of $PGL_{n+1}{\Bbb R}$,
we assume $Z_1$ and  $Z_2$ intersect transversally.
Components of $Z_1$ and $Z_2$ are Schubert varieties defined
by real flags.
Moreover, each component of $Z_1$ has complementary dimension to
each component of  $Z_2$.
In a  flag variety,
Schubert varieties of complimentary dimension
which meet transversally and are defined by real flags
either have empty intersection, or meet in
a single real point.
Thus $Z_1\bigcap Z_2$ consists entirely of real points.

By Lemma~\ref{lemma:curve_selection}, for each $i=1,2$, there is a real
algebraic map $f_i : C_i \rightarrow \overline{\phi(V_i)}$
where $C_i$ is a smooth curve, and  a point $s_i$ on a connected arc $S_i$
of $C_i({\Bbb R})$ such that $f(S_i -\{s_i\})\subset\phi(V_i({\Bbb R}))$
and $f_i(s_i) = Z_i$.

The enumerative problem $\Psi'\rightarrow V'$ given by
$f_1^*\Psi_1\rightarrow C_1$ and $f_2^*\Psi_2\rightarrow C_2$
is fully real, as $\Psi'_{(s_1,s_2)} = Z_1\bigcap Z_2$.
Since $(s_1,s_2)\in
\overline{(S_1-\{s_1\})\times(S_2-\{s_2\})\bigcap V'({\Bbb R})}$,
there is a point $(t_1,t_2)\in (S_1-\{s_1\})\times(S_2-\{s_2\})$
such that
$\Psi'_{(t_1,t_2)} = (f_1^*\Psi_1)_{t_1}\bigcap(f_2^*\Psi_2)_{t_2}$
is transverse and consists entirely of real points.
Since   $f_i(t_i)\in \phi(V_i({\Bbb R}))$,
we see that  $(f_i^*\Psi_i)_{t_i} = (\Psi_i)_{f_i(t_i)}$ is a fibre
of $\Psi_i$ over a point of $V_i({\Bbb R})$.
This shows $\Psi\rightarrow V$ is fully real.
\QED

\section{Powers of Enumerative Problems}

A method to construct a new fully real enumerative problem out of
a given one is illustrated by a proof of B\'ezout's Theorem
in the plane.
We will formalize this method.

\subsection{B\'ezout's Theorem}
{\em
Let $d_1$ and $d_2$ be positive integers.
Then there exist smooth real plane curves $D_1$ and $D_2$
of degrees $d_1$ and $d_2$  meeting transversally in
\medskip$d_1\cdot d_2$ real points.}

\noindent{\bf Proof:}
Two distinct real lines meet in a single real point.
Thus if $D_1$
consists of $d_1$ distinct real lines, $D_2$ of $d_2$,
and if $D_1$ and $D_2$ meet transversally, then $D_1\bigcap D_2$ is
$d_1\cdot d_2$ real points.

The family of real reduced degree $d$ plane
curves has general member a smooth curve and contains
all cycles of $d$ distinct real lines.
Thus the enumerative problem of intersecting reduced curves $D_1$
and $D_2$ of respective degrees $d_1$ and $d_2$ is fully real of
degree $d_1\cdot d_2$.
Moreover, pairs of smooth real curves are dense in the set of pairs
of reduced real curves, showing the enumerative problem of intersecting
two smooth plane curves of respective degrees  $d_1$ and $d_2$ is
fully real of
degree $d_1\cdot d_2$.
\QED

\subsection{Powers of intersection problems}
Suppose $\Xi\rightarrow U$ is a family of multiplicity
free cycles on $X$ and $d$ is a positive integer.
If the locus of $d$-tuples
$(u_1,\ldots,u_d)$ such that no two of $\Xi_{u_1},\ldots,\Xi_{u_d}$
share a component is dense in $U^d$, then let
$U^{(d)}$ be an open subset of that locus.
Let $\Xi^{\oplus d}\rightarrow U^{(d)}$ be the family
of multiplicity free cycles whose fibre over
$(u_1,\ldots,u_d)\in U^{(d)}$ is $\sum_{j=1}^d \Xi_{u_j}$.

Suppose $\Xi_1\rightarrow U_1,\ldots,\Xi_b\rightarrow U_b$ are families
of multiplicity free cycles on $X$ giving an
intersection problem $\Xi\rightarrow U$
and $d_1,\ldots,d_b$ is a sequence of positive integers.
Then the families $\Xi_1^{\oplus d_1}\rightarrow U_{1}^{(d_1)},\ldots,
\Xi_b^{\oplus d_b}\rightarrow U_{s}^{(d_b)}$
give a well-posed intersection problem if
general members of the  famililies
$\Xi\rightarrow U$ and $\Xi_i\rightarrow U_i$
meet properly, for $1\leq i\leq b$.

When a reductive group $G$ acts transitively
on  $X$ and the families of cycles are $G$-stable,
$\Xi_1^{\oplus d_i},\ldots,\Xi_b^{\oplus d_b}$ give an intersection
problem.
Moreover, if $\Xi\rightarrow U$ is fully real, then
so is that intersection problem.
We produce a witness with a particular form.

\subsection{Lemma}\label{lemma:transverse_technicality}
{\em
Suppose $\Xi_1\rightarrow U_1,\ldots,\Xi_b\rightarrow U_b$
give a fully real enumerative problem of degree $d$.
Let $d_1,\ldots,d_b$ be a sequence of positive integers
and suppose that for $1\leq i\leq b$, $V_i$ is  $G$-stable subset of
$U_{i}^{(d_i)}$ such that
$\Delta^{d_i}U_i({\Bbb R})\subset \overline{V_i({\Bbb R})}$,
as subsets of $U_i({\Bbb R})^{d_i}$.

Then for $1\leq i\leq b$, there exists $v_i\in V_i({\Bbb R})$ such that
$(\Xi_1^{\oplus d_1})_{v_1},\ldots,(\Xi_b^{\oplus d_b})_{v_b}$
intersect transversally in
$d\cdot d_1\cdots d_b$ real points.}
\medskip

\noindent{\bf Proof:}
The restriction $\Psi_i$ of $\Xi_i^{\oplus d_i}$ to $V_i$ is
$G$-stable.
Thus $\Psi_1,\ldots,\Psi_b$ give a well-posed
enumerative problem $\Psi\rightarrow V$.
We show this is  fully real and compute its degree.

Since $\Xi\rightarrow U$ is fully real, there
is an open subset $R$ of points
$u \in U({\Bbb R})$ such that $\Xi_u$ is
$d$ distinct real points.
Since $U({\Bbb R})\subset \prod_{i=1}^b U_i({\Bbb R})$,
for $1\leq i\leq b$ there exists an open subset
$R_i$ of $U_i({\Bbb R})$ such that $\prod_{i=1}^b R_i \subset R$.
Then
$V_i({\Bbb R})\bigcap R^{d_i}\neq \emptyset$,  as
$\Delta^{d_i} R_i\subset \Delta^{d_i}U_i({\Bbb R})
\subset \overline{V_i({\Bbb R})}$.
Thus $R' := V({\Bbb R})
\bigcap \prod_{i=1}^b R_i^{d_i}$
is nonempty,
as $V({\Bbb R})$ is  dense in
$ \prod_{i=1}^b V_i({\Bbb R})$.

Let $w =(w_{11},\ldots,w_{1d_1},\ldots,w_{b1}\ldots,w_{bd_b})\in R'$.
Here, $w_{ij}\in R_i$ and
$(w_{i1},\ldots,w_{id_i})\in  V_i({\Bbb R})$.
If $1\leq j_i\leq d_i$,  
then $(w_{1j_1},\ldots,w_{bj_b})\in U({\Bbb R})$.
Furthermore, $\Psi_w = \bigcap_{i=1}^b (\Psi_i)_{(w_{i1},\ldots,w_{id_i})}$
is a transverse intersection, as $R' \subset V$.
Since $(\Psi_i)_{(w_{i1},\ldots,w_{id_i})} =
\sum_{j=1}^{d_i} (\Xi_i)_{w_{ij}}$, we have
$$
\Psi_w \ =\ \bigcap_{i=1}^b \,\sum_{j=1}^{d_i} (\Xi_i)_{w_{ij}}
\ =\ \sum_{\stackrel{\mbox{\scriptsize$j_1,\ldots,j_b$}}{1\leq j_i\leq d_i}}
\ \bigcap_{i=1}^b\,(\Xi_i)_{w_{ij_i}}
\ =\ \sum_{\stackrel{\mbox{\scriptsize$j_1,\ldots,j_b$}}{1\leq j_i\leq d_i}}
\Xi_{(w_{1j_1},\ldots,w_{bj_b})}.
$$
Since this  intersection is transverse,
it consists of  $d\cdot d_1\cdots d_b$ real points.
\QED

\subsection{Real B\'ezout's Theorem} {\em
Let $d_1,\ldots,d_b$ be positive integers.
Then there exist smooth real hypersurfaces
$H_1,\ldots,H_b$ in ${\Bbb P}^b$ of respective degrees
$d_1,\ldots,d_b$ which intersect transversally
in  $d_1\cdots d_b$ real points.
}\medskip

\noindent{\bf Proof:}
Let $\Xi\rightarrow U$ be the family of hyperplanes in ${\Bbb P}^b$.
Since $b$ real hyperplanes in general position meet in a real point,
either simple checking or Lemma~\ref{lemma:transverse_technicality} with
$V:=U^{(d_i)}$ shows that $\Xi^{\oplus d_1},\ldots,\Xi^{\oplus d_b}$
give a fully real enumerative problem of degree $d_1\cdots d_b$.
Note that $\Xi^{\oplus d_i}\rightarrow U^{(d_i)}$ is the family of
hypersurfaces composed of  $d_i$ distinct hyperplanes.

Let $W_i\subset {\Bbb P}(\mbox{\em Sym}^{d_i}{\Bbb C}\,^{b+1})$ be the space
of forms of degree $d_i$ with no repeated factors and
$\Gamma_i\rightarrow W_i$ the family of reduced degree $d_i$ hypersurfaces.
Let $W'_i\subset W_i$ be the dense subset of forms determining
smooth hypersurfaces.
Note that $U^{(d_i)}\subset W_i$ and
$\Xi^{\oplus d_i} = \Gamma_i|_{U^{(d_i)}}$.

It follows that $\Gamma_1,\ldots,\Gamma_b$ give a fully real
enumerative problem of degree $d_1\cdots d_b$.
Let $R$ be an open set of witnesses.
Since $U^{(d_i)}({\Bbb R})\subset\overline{W'_i({\Bbb R})}$
and $R$ meets $\prod_{i=1}^sU^{(d_i)}({\Bbb R})$,  we see that
$R$ meets $\prod_{i=1}^sW'_i({\Bbb R})$.
That is,
there exist smooth real hypersurfaces
$H_1,\ldots,H_b$ in ${\Bbb P}^b$ of respective degrees
$d_1,\ldots,d_b$ which intersect transversally
in  $d_1\cdots d_b$ real points.
\QED

\section{$(n-2)$-planes meeting rational normal curves in
${\Bbb P}^n$}

Let ${\Bbb G}_{n-2}{\Bbb P}^n$ be the Grassmannian of $(n-2)$-planes in
${\Bbb P}^n$, a variety of dimension $2n-2$.
Those  $(n-2)$-planes which meet a curve form a hypersurface
in ${\Bbb G}_{n-2}{\Bbb P}^n$.
We synthesize ideas of previous sections to
prove the following theorem.

\subsection{Theorem}\label{thm:rational_normal}
{\em
The enumerative problem of $\,(n-2)$-planes meeting $2n-2$ general
rational normal
curves in ${\Bbb P}^n$ is fully real  and has degree
${2n-2\choose n-1}n^{2n-3}$.
}\medskip

\noindent{\bf Proof:}
Identifying ${\Bbb P}^n$ with its dual projective space
gives an isomorphism
${\Bbb G}_{n-2}{\Bbb P}^n
\stackrel{\sim}{\rightarrow} {\Bbb G}_1{\Bbb P}^n$,
mapping Schubert subvarieties to Schubert subvarieties.
By Theorem~C of~\cite{sottile_real_lines}, any enumerative problem
involving Schubert subvarieties of ${\Bbb G}_{n-2}{\Bbb P}^n$ is fully real.
In particular, the enumerative problem given by $2n-2$ copies of the family
$\Xi \rightarrow U$ is fully real,
where $U = {\Bbb G}_1{\Bbb P}^n$ and the fibre of $\Xi$
over $l\in U$ is the Schubert variety $\Omega_l$
of $(n-2)$-planes meeting $l$.

We compute its degree, $d$.
The image of $\Omega_l$
under the isomorphism ${\Bbb G}_{n-2}{\Bbb P}^n
\stackrel{\sim}{\rightarrow} {\Bbb G}_1{\Bbb P}^n$
is the Schubert subvariety of all lines meeting
a fixed $(n-2)$-plane.
Thus $d$ is the number of lines meeting $2n-2$ general $(n-2)$-planes
in ${\Bbb P}^n$.
By Corollary 3.3 of~\cite{sottile_real_lines}, this is the number of
(standard) Young tableaux of shape $(n-1,n-1)$, which is
$\frac{1}{n} {2n-2\choose n-1}$, by the hook length formula
of Frame, Robinson, and Thrall~\cite{FRT}.

Let $e_0,\ldots,e_n$ be real points spanning ${\Bbb P}^n$.
For $1\leq i\leq n$, let $l_i := \Span{e_{i-1},e_i}$.
Then $\Omega_{l_1}+ \cdots+\Omega_{l_n}$
is the fibre of $\Xi^{\oplus n}$ over $(l_1,\ldots,l_n)\in U^{(n)}({\Bbb R})$.
Let $V = PGL_{n+1}{\Bbb C}\cdot (l_1,\ldots,l_n)\subset U^{(n)}$.
For $t\in [0,1]$ and $1\leq i\leq n$, let
$$
l_i(t) \ :=\ \Span{te_{i-1} + (1-t)e_{\overline{i -1}},\,
t e_i + (1-t) e_{\overline{i}}},
$$
where $\overline{j}\in \{0,1\}$ is congruent to $j$ modulo 2.
Let $\gamma(t) := (l_1(t),\ldots,l_n(t))$.
If $t\in (0,1]$, then $\gamma(t)\in V({\Bbb R})$.
Since $\gamma(0) = (l_1,\ldots,l_1)$ and
$\Delta^nU({\Bbb R}) = PGL_{n+1}{\Bbb R}\cdot \gamma(0)$,
it follows that $\Delta^nU({\Bbb R})\subset \overline{V{(\Bbb R})}$.
Then, by Lemma~\ref{lemma:transverse_technicality}, there exist points
$v_1,\ldots,v_{2n-2} \in V({\Bbb R})$ such that
$\Xi^{\oplus n}_{v_1},\ldots,\Xi^{\oplus n}_{v_{2n-2}}$
meet transversally in
${2n-2\choose n-1}n^{2n-3}$ points.

Let $p(m):= n\cdot m +1$, the Hilbert polynomial of a rational normal curve
in ${\Bbb P}^n$.
Let ${\cal H}$ be the open subset of the Hilbert scheme
parameterizing reduced schemes with Hilbert polynomial $p$.
Let $\Psi\subset {\cal H}\times {\Bbb G}_{n-2}{\Bbb P}^n$
be the family of multiplicity free cycles on
${\Bbb G}_{n-2}{\Bbb P}^n$
whose fibre over a curve $C\in {\cal H}$ is the hypersurface
of $(n-2)$-planes meeting $C$.

Note that $p$ is also the Hilbert polynomial of
$l_1\bigcup\cdots\bigcup l_n$.
Let $\lambda \in {\cal H}$ be the point representing
$l_1\bigcup\cdots\bigcup l_n$.
If $V'$ is the $PGL_{n+1}{\Bbb C}$-orbit of $\lambda$ in ${\cal H}$,
then $\Psi|_{V'}\rightarrow V'$ is isomorphic to the family
$\Xi^{\oplus n}\rightarrow V$, under the obvious isomorphism
between $V$ and $V'$.
It follows that the enumerative problem given by $2n-2$
copies of $\Psi\rightarrow {\cal H}$ is fully real.

Let $W$ be the subset of ${\cal H}$ representing rational
normal curves.
We claim $V'({\Bbb R})\subset \overline{W({\Bbb R})}$, from which
it follows that
the enumerative problem of $(n-2)$-planes meeting $2n-2$ rational
normal curves in ${\Bbb P}^n$ is fully real and has degree
${2n-2\choose n-1}n^{2n-3}$.

Let $[x_0,\ldots,x_n]$ be homogeneous coordinates for ${\Bbb P}^n$
dual to the basis $e_0,\ldots,e_n$.
For $t\in {\Bbb C}$, define the ideal ${\cal I}_t$ by
$$
{\cal I}_t \ :=\
( x_ix_j - t x_{i+1}x_{j-1}\, |\, 0\leq i<j\leq n\ \mbox{and}\ j-i\geq 2).
$$
For $t\neq 0$, ${\cal I}_t$ is the ideal of a rational normal curve
and ${\cal I}_0$ is the ideal of $l_1\bigcup\cdots\bigcup l_n$.

This family of ideals is flat.
Let $\varphi:{\Bbb C}\rightarrow {\cal H}$
be the map representing this family.
Then $\varphi({\Bbb R}-\{0\})\subset W({\Bbb R})$.
Noting $\varphi(0)=\lambda$ shows
$\lambda \in \overline{W({\Bbb R})}$.
Since $W({\Bbb R})$ is $PGL_{n+1}{\Bbb R}$-stable,
we conclude that
$V'({\Bbb R})\subset \overline{W({\Bbb R})}$.
\QED

\end{document}